\documentclass{amsart}
\usepackage{xspace}
\usepackage[dvipsnames]{xcolor}
\usepackage{amsmath}
\usepackage{amssymb}
\usepackage{enumitem}

\usepackage{graphicx}
\usepackage{epstopdf}
\usepackage{tikz}
\usepackage{subfigure}
\usepackage{siunitx}
\DeclareSIUnit\darcy{d}
\usetikzlibrary{positioning}

\usepackage{bm}
\usepackage{enumitem}
\setlist{nolistsep}

\newcommand{\ld}{\ell}
\newcommand{\es}{e}
\newcommand{\ra}{r}

\newcommand{\trans}{T}
\newcommand{\bigand}{\quad\text{and}\quad}
\newcommand{\epsi}{\ensuremath{\varepsilon}}

\newcommand{\Real}{\mathbb{R}}
\newcommand{\Pcal}{\mathcal{P}}
\newcommand{\Kcal}{\mathcal{K}}
\newcommand{\grad}{\nabla}

\newcommand{\dive}{\grad\cdot }

\renewcommand{\vec}[1]{\bm{#1}}
\newcommand{\abs}[1]{\left|#1\right|}
\newcommand{\ddiv}[1]{\text{div}\ifx\empty #1\else _{\hspace{-1pt}#1}\xspace}
\newcommand{\dgrad}[1]{\text{grad}\ifx\empty #1\else _{\hspace{-1pt}#1}\xspace}

\newcommand{\fracpar}[2]{\frac{\partial #1}{\partial #2}}
\DeclareMathOperator{\id}{I}

\DeclareMathOperator{\trace}{tr}

\usepackage{soul}
\soulregister\cite7
\soulregister\ref7
\soulregister\eqref7
\soulregister\num7
\soulregister\pageref7

\author{Halvor M{\o}ll Nilsen \and Idar Larsen \and Xavier Raynaud} \title[Combining MDEM and VEM for fracturing modeling]{Combining the Modified Discrete Element Method
  with the Virtual Element Method for Fracturing of Porous Media.}

\begin{document}
\maketitle

\begin{abstract}
  Simulation of fracturing processes in porous rocks can be divided into two main
  branches: (i) modeling the rock as a continuum which is enhanced with special
  features to account for fractures, or (ii) modeling the rock by a discrete (or
  discontinuous) approach that describes the material directly as a collection of
  separate blocks or particles, e.g., as in the discrete element method (DEM).  In
  the modified discrete element (MDEM) method, the effective forces between virtual
  particles are modified in all regions, without failing elements, so that they
  reproduce the discretization of a first order finite element method (FEM) for
  linear elasticity. This provides an expression of the virtual forces in terms of
  general Hook's macro-parameters. Previously, MDEM has been formulated through an
  analogy with linear elements for FEM. We show the connection between MDEM and the
  virtual element method (VEM), which is a generalization of FEM to polyhedral
  grids. Unlike standard FEM, which computes strain-states in a reference space, MDEM
  and VEM compute stress-states directly in real space.  This connection leads us to
  a new derivation of the MDEM method. Moreover, it gives the basis for coupling
  (M)DEM to domain with linear elasticity described by polyhedral grids, which makes
  it easier to apply realistic boundary conditions in hydraulic-fracturing
  simulations. This approach also makes it possible to combine fine-scale (M)DEM
  behavior near the fracturing region with linear elasticity on complex reservoir
  grids in the far-field region without regridding. To demonstrate the simulation of
  hydraulic fracturing, the coupled (M)DEM-VEM method is implemented using the Matlab
  Reservoir Simulation Toolbox (MRST) and linked to an industry-standard reservoir
  simulator. Similar approaches have been presented previously using standard FEM,
  but due to the similarities in the approaches of VEM and MDEM, our work is a more
  uniform approach and extends these previous works to general polyhedral grids for
  the non-fracturing domain.
\end{abstract}

\clearpage

\section{Introduction}

Effective control of flows in geological formations is a key factor for exploiting
resources that are highly important for the society, such as ground water, geothermal
energy, geological storage of CO2, high quality fossil fuel (gas and oil) and
potentially large scale storage of energy in terms or heat or gas. Today 60\% of the
world energy consumption is based on oil and natural gas resources \cite{iea2015}. In
addition 19\% is based on coal which needs large scale CO2 storage to be safely
exploited without large scale impact on climate \cite{ipcc2014}. Geothermal energy is
an important source of green energy, which would be even more valuable in the future
as the supply of fossil fuel is expected to decrease. Gas storage is today an
integrated part of the energy supply and provides reliable large scale storage of
energy. It enables to both attenuate the volatility of energy prices and ensure
energy security.

The use of all of the above resources will benefit from a reliable control of the
flow properties around the wells that are used to exploit them. Increased injectivity
is particularly important for exploiting resources in tight formation or where high
flow rates are required.  For enhanced geothermal applications, rock fracturing is a
prerequisite for economical exploitation. For CO2 injection, where large volumes of
fluid have to be injected, high injectivity limits the increase in pressure near the
well and simplifies the operation. Tight formations contain much of the hydrocarbon
reserves. The exploitation of these formations has been a driving force for the
technology of fracking, which is a more drastic well stimulation than the traditional
ones. When increasing the injectivity in a well, it is of vital importance to be able
to predict and control fracturing to avoid unwanted fractures or even induced seismic
events, which may cause environmental damage as well as the disruption of
operations. The key for enabling high injectivity is to induce and control the
fracturing process using the coupling between fluid flow, heat and rock
mechanics. The failure of the rock and the propagation of fractures depend on both
global and local effects, through the stress distribution which is intrinsically
global and failure criterias which are local. It is therefore important to have
flexible simulation tools that are able to cover both large scale features with
complex geometry and include specific fracture dynamics where the fracturing
processes occur.

Numerical methods for simulating fracturing can typically be classified either as
continuum or discontinuum based methods \cite{Jing2002409}. The modeling of
fracturing in brittle materials like rock is particularly difficult. It is determined
by the stress field in the vicinity of the fracture tip. As \cite{Griffith163}
showed, failures happen when the global energy release is larger than the energy
required by the fracturing process. The first depends on the global stress field
while the latter is associated with the energy needed at small scale to create a
fracture. For brittle materials where the failure happens at very small scales,
linear elasticity governs the behavior in most of the domain but the solution of
linear elasticity in the presence of fractures is singular near the fracture tip (see
\cite{Kuna2013} for a general description). This introduces challenges for numerical
calculations and often results in artificial grid dependence of the simulated
dynamics. From a physical point of view, such effects are removed when plasticity is
introduced, however the length scale of this region may be prohibitively small to be
resolved numerically on the original model.  Several techniques have been introduced
to incorporate the singularity at the fracture tip into the numerical calculations
explicitly, for example specific tip elements in the finite element FEM method. In
general, the methods using global energy arguments are less sensitive to the choice
of numerical methods than those that use estimates of the strength of the singularity
\cite{Kuna2013}. In the case of the fracturing of natural rock, the uncertainty in
the model is large, small scale heterogeneities are important and several different
fracturing mechanisms complicate the structure. Discrete modeling techniques have
been very successful in this area \cite{Lisjak2014301}, in particular if complex
behaviors should be simulated.

An essential component to model fracturing is therefore the ability to account
in a flexible manner for both large and small scale behaviors. This is reflected
in the widespread use of tools based on analytic models for hydraulic fracturing
(for a review see \cite{Li20158}). However, it becomes a challenge to
incorporate interaction of fractures and fine-scale features into those
simulation tools. Because of their simplicity and flexibility for incorporating
different fracturing mechanisms, discrete element methods (DEM), also called in
their explicit variants, distinct element methods, have been one of the main
techniques used for hydraulic fracturing in commercial
simulators. 
Those methods exploit the ability of easily modifying the interactions and
connections between the discrete particles or elements. For continuum models, such
behavior is more difficult to account for. However, the parameters in the DEM model
are not directly related to physical macro-scale parameters and there are
restrictions on the range of the parameters that can be simulated. In particular,
\cite{Alassi2012} showed that only Poisson's ratios (in plain stress) smaller than
$1/3$ can be considered. The modified discrete element method (MDEM) was introduced
to get rid of this restriction and gives also easier relationships between the macro
parameters in the linear elastic domain, while keeping the advantages of DEM in the
treatment of fracture. In this work, we show the connection between MDEM and the
recent development of Virtual Element Methods (VEM). Such approach provides a simple
derivation of the MDEM framework, and also highlight the discrepancy of the original
DEM model from linear elasticity. The linear version of the VEM methods for
elasticity can be used to extend first-order FEM on simplex grids, which was the
basis of the MDEM method, to general polyhedral grids. We use the fact that both DEM,
MDEM and VEM share the same degrees of freedom in the case of simplex grids to derive
smooth couplings between these methods. A similar approach has been followed for
coupling FEM with DEM previously \cite{Pan1991,Billaux2004}. The introduction of VEM
opens the possibility for flexible gridding on general polyhedral grids in the
far-field region while keeping the DEM/MDEM flexibility in the near fracture
domains. Geological formations are typically the result of deposition and erosion
processes, which lead to layered structures and faults. Geometrical models using
polyhedral grids, such as Corner Point Grid \cite{ponting1989corner}, Skua Grid
\cite{Gringarten2008} and Cut-Cell \cite{Mallison2014} are natural in this context
and correspond to grids used in the industry of flow modeling in reservoirs. Our
proposed method therefore may simplify the incorporation of fracture simulation in
realistic subsurface applications.

\section{Methods}
We study the methods for the standard equations of linear elasticity given by
\begin{equation}
  \label{eq:lin_elast_cont}
  \begin{aligned}
    \dive \sigma &= \vec{f},\\
    \epsi &= \frac{1}{2}(\grad +\grad^{T})\vec{ u},\\
    \sigma &= C \epsi,
  \end{aligned}
\end{equation}
where $\sigma$ is the Cauchy stress tensor, $\epsi$ the infinitesimal strain tensor
and $u$ the displacement field. The linear operator $C$ is the fourth-order Cauchy
stiffness tensor. In Kelvin notation \cite{GPR:GPR1049}, a three-dimensional
symmetric tensor $\{\epsi_{ij}\}$ is represented as an element of $\Real^6$ with
components
\begin{equation*}
  [\epsi_{11},\epsi_{22},\epsi_{33}, \sqrt{2} \epsi_{23}, \sqrt{2}\epsi_{13}, \sqrt{2} \epsi_{12}]^{\trans}
\end{equation*}
while a two-dimensional symmetric
tensor is represented by a vector in $\Real^3$ given by $[\epsi_{11},\epsi_{22},
\sqrt{2} \epsi_{12}]^{\trans}$. Using this notation $C$ can be represented by a
$6\times 6$ matrix and inner products of tensors correspond to the normal
inner-product of vectors. For isotropic materials, we have the constitutive equations
\begin{equation}
  \label{eq:isosigma}
  \sigma = 2\mu\epsi + \lambda\trace(\epsi)\id.
\end{equation}
where $\mu$ and $\lambda$ denote the Lam\'e constants. The elastic energy density is
given by $\frac12\sigma : \epsi$ where we use the standard scalar product for
matrices defined as
\begin{equation*}
  \alpha:\beta = \trace(\alpha^{\trans}\beta) = \sum_{i,j=1}^{3}\alpha_{i,j}\beta_{i,j},
\end{equation*}
for any two matrices $\alpha,\beta \in \Real^{3\times 3}$.

\subsection{Discete element method}

Discrete element methods consist of modeling the mechanical behaviour of a continuum
material by representing it as a set of particles, or discrete elements. The forces
in the material are then modeled as interaction forces between the particles. There
are several variants of the discrete element method \cite{Lisjak2014301}.  Here we
will use the simple version introduced in \cite{Alassi2012} where the particles are
discs in 2D and spheres in 3D. We will also restrict the treatment to the linear case
to compare with linear elasticity, but this is not a restriction of the
method. 
For more in-depth presentation of different variants see \cite{Lisjak2014301} or
\cite{Cundall1979} and the references therein.

The starting point of the DEM methods has its background in the description of
granular media. This has a long history starting from the description of the contact
force by Hertz and \cite{Mindlin1949}.  In this field an important question was to
study how the effective elastic modulus of the bulk was related to the microscopic
description \cite{Dvorkin1996,Walton1987}. In DEM the basic ideas is to use a
microscopic description to simulate the behavior of the bulk modulus.
\begin{figure}
\begin{center}
\subfigure[DEM]{
\includegraphics{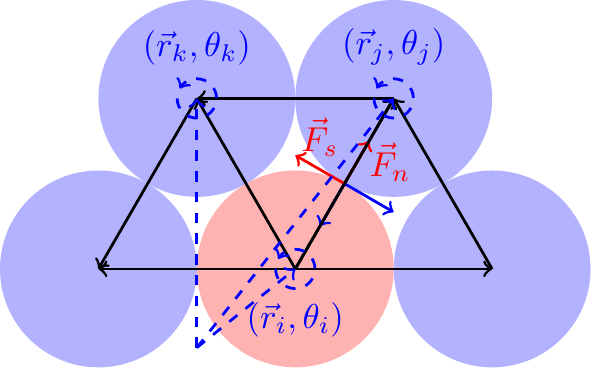}\label{fig:dem_mdem1}
}
\subfigure[MDEM]{
\includegraphics{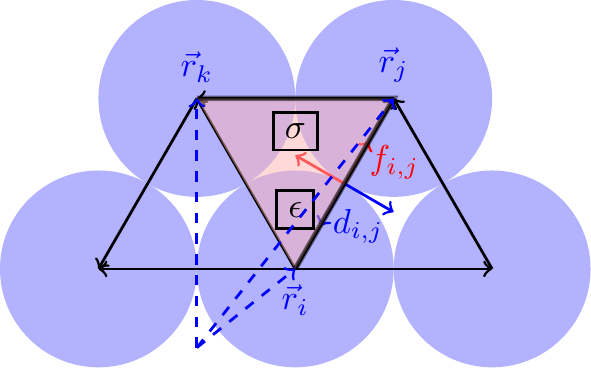}\label{fig:dem_mdem2}
}

\end{center}
\caption{This show the fundamental quantities used in DEM \subref{fig:dem_mdem1} and
  MDEM \subref{fig:dem_mdem2}. The unit colored in red for MDEM defines the area
  which is used to calculated the one sided forces.}\label{fig:dem_mdem}
\end{figure}

For a complete relation between general DEM method and linear elasticity using shear
forces, it is necessary to introduce the material laws for micropolar materials, see
\cite{Alassi2012}. This introduces an extra variable associated with local rotation,
as illustrated in Figure \ref{fig:dem_mdem}. For an isotropic
micropolar 
media the strain stress relation is
\begin{equation}
  \label{eq:defstressmicro}
  \sigma = 2\mu\epsi + \lambda\trace(\epsi)\id + \kappa   (\tau- \phi)
\end{equation}
where the extra variable $\phi$ describe the local rotation and $\tau$ represent
the asymmetric part of the strain tensor, i.e, rigid body rotations. In terms of
displacements it can be written
\begin{equation}
  \tau = \frac{1}{2}(\grad u - \grad^{\trans} u).
\end{equation}
The state variable are the displacement $u$ and the local rotation $\phi$ and the
total elastic energy is given by
\begin{equation}
  \label{eq:deftotE}
  E =\int_\Omega(\mu\epsi:\epsi + \frac{\lambda}2\trace(\epsi)^2 + \frac{\kappa}{2}(\tau- \phi):(\tau- \phi))\,dx.
\end{equation}
By computing the variation of $E$,
\begin{equation*}
  \delta E =\int_\Omega(\dive(2\mu\epsi + \lambda\trace(\epsi)\id + \kappa(\tau-
  \phi))\cdot\delta u)\,dx
  + \int_\Omega \kappa(\tau - \phi):\delta\phi\,dx,
\end{equation*}
we obtain the governing equations of the system, that is, the linear momentum conservation equation,
\begin{subequations}
  \label{eq:goveqmicro}
  \begin{equation}
    \label{eq:linmomconsmicro}
    \dive(2\mu\epsi + \lambda\trace(\epsi)\id + \kappa(\tau-
    \phi)) = 0
  \end{equation}
  and the angular momentum conservation equation
  \begin{equation}
    \label{eq:angmomconsmicro}
    \tau- \phi = 0.
  \end{equation}
\end{subequations}
We will use the expression of the stress given by \eqref{eq:defstressmicro} to
compare with the DEM model which include local rotations. Let us consider two
particles $p_1$ and $p_2$ which are connected through a contact denoted $m$. For a
particle $p_i$, $i=1,2$, we denote by $X^{p_i}$ the position of the particle. The
degrees of freedom of the system are the displacement $U^{p_i}$ and the microrotation
$\theta^{p_i}$ for each particle $p_i$. Furthermore we set
\begin{equation*}
  \Delta X^m = X^{p_2} - X^{p_1}, \quad \Delta U^m = U^{p_2} - U^{p_1},\quad I^m = \frac{\Delta X^m}{\abs{\Delta X^m}}.
\end{equation*}
We introduce also the distance between the particles, $d_m = \abs{\Delta X^m}$. We
use the cross-product to represent the action of a rotation so that the rotation
given by a vector $\theta$ is the mapping given by $X\to\theta\times X$. Our
description of DEM follows \cite{Alassi2012} with slight differences in the
notation. We introduce the normal and shear forces,
\begin{equation}
  F_{n}^m = k_n \Delta U_{n}^m \quad F_{s}^m=k_s \Delta U_{s}^m.
\end{equation}
For a given contact $m$ the relative shear and normal displacement are given by
\begin{align*}
  \Delta U_{n}^m &= (\Delta U^m\cdot I^m)I^m,\\
  \Delta U_{s}^m &= \Delta U^m - \Delta U_{n}^m - \theta^{m}\times\Delta X^m,
\end{align*}
where $\theta^m = \frac12(\theta^{p_1} + \theta^{p_2})$. Note that in the case where
two adjacent spheres roll one over the other without sliding, we have
$\theta^{p_1}=-\theta^{p_2}$ so that the term $\Delta U_s^m$ accounts only for the
sliding part of the tangential component. Let us define the total force over a
contact $m$ as
\begin{equation*}
  F^m =F_{n}^m +  F_{s}^m.
\end{equation*}
Using the definition of stress tensor $\sigma$, we have that, at the contact between the
spheres of center $p_2$ and $p_1$, and assuming that there exists a non-zero contact
surface $ds$, the force $F^m$ can be written as
\begin{equation*}
  F^m = \sigma I^m ds,
\end{equation*}
as $I^m$ points in the normal direction. The Cauchy's formula for the stress matrix
\cite{Alassi2012}, which is meant to invert this relation, is given by
\begin{equation*}
  \sigma = \frac{1}{V}\sum_{m = 1}^{N_c} d_m F^m\otimes I^m,
\end{equation*}
where $N_c$ denotes the number of contact points, that is the number of spheres in
contact. Let us consider a linear deformation $\ld$ and write $\Delta U^m$ as
\begin{equation}
  \label{eq:lindisp}
  \Delta U^m = \ld \Delta X^m = (\es +  \ra)\Delta X^m,
\end{equation}
where the tensor $\es$ and $r$ are respectively the symmetric and skew-symmetric
parts of $\ld$. Since $\ra$ is a skew-symmetric matrix, it corresponds to a rotation
and, abusing the notations, we will write indifferently $\ra\Delta X^m$ or $\ra\times
\Delta X^m$ to denote the same rotation operator (here applied to $\Delta X^m$). To
proceed with the identification of the stress tensor, we assume small displacement,
that is, $e$ and $r$ are small compared with the identity, and we assume also
$\theta^m=\theta$ for some constant $\theta$. We use \eqref{eq:lindisp} and obtain
\begin{equation}
  \label{eq:Fnm}
  F_n^m = d_mk_n(I^{m}\cdot \es I^m + I^{m}\cdot \ra I^m) I^m = d_mk_n(I^m\cdot \es I^m) I^m,
\end{equation}
as $\ra$ is skew-symmetric. For $F_s^m$, we have
\begin{equation}
  \label{eq:Fns}
  F_s^m = d_mk_s( \es I^m - I^m\cdot \es I^m + (\ra - \theta)\times I^m ).
\end{equation}
Hence, we obtain the following expression for the stress tensor,
\begin{multline}
  \label{eq:sigma_dem}
  \sigma = \frac{1}{V}\sum_{m = 1}^{N_c}  d_m^2\big((k_n - k_s)(I^m\cdot \es I^m)\,[I^m\otimes I^m] \\
    + k_s\,\left([(\es I^m)\otimes I^m] + [((\ra - \theta)\times I^m) \otimes
      I^m]\right)\big)
\end{multline}
To illustrate the restriction that this expression imposes on the parameters, we
consider a square packing in 3D. In this case, there are $N_c = 6$ contact points
and, using $I^1 = -I^2 = (1,0,0)$, $I^3 = -I^4 = (0,1,0)$ and $I^5 = -I^6 = (0,0,1)$,
we obtain
\begin{equation*}
  \sigma = 2(k_n - k_s)\trace(\es)\id + 2k_s\es + 2 k_s (\ra - \theta).
\end{equation*}
Note that we do not take $V$ equal to the volume of the sphere but $V=1$, that is the
effective volume. In the expression above, $\ra$ and $\theta$ must be seen as
matrices and not as vectors as in \eqref{eq:Fns}. We can also identified the
parameter $\phi$ of local rotation for micropolar media with the local rotation
$\theta$ in the DEM model. This gives the following Lam\'e coefficients,
\begin{equation*}
  \lambda = 2(k_n - k_s) \bigand \mu = k_s.
\end{equation*}
Hence, as $\sigma:\es=2\mu\sum_{ij}\es_{ij}^2 + \lambda\sum_{i}\es_{ii}^2 =
2k_s\sum_{i\neq j}\es_{ij}^2 + 2k_n\sum_{i}\es_{ii}^2 $, we can conclude that, for
square lattices, this is only stable if $k_s>0$ (we assume $k_n,k_s\geq 0$). However,
this is not a restriction for simplex grids.  Using the same approach as above but
now for regular simplices, it is shown in \cite{Alassi2012} that
\begin{equation}
 \mu= k_n + k_s \quad
 \lambda =  k_n - k_s
\end{equation}
Since $k_s$ and $k_n$ are naturally positive, this restricts the Poisson's ratio to
\begin{equation}
  \nu = \frac{\lambda}{2(\mu + \lambda)}  = \frac{1}{4} (1-\frac{k_s}{k_n})  < \frac{1}{4}
\end{equation}
in the 3D case. For the 2D case, we obtain the same expression in the case of plane
strain boundary conditions and, in the case of plane stress, we get
\begin{equation}
 \nu = \frac{\lambda}{2 \mu + \lambda}  = \frac{k_n - k_s}{3 k_n + k_s} = \frac{1}{3} \frac{1-\frac{k_s}{k_n}}{1+\frac{k_s}{k_n}},
\end{equation}
which implies $-1<\nu<\frac13$. These limitations on the physical parameters have
been the main motivation for introducing MDEM. Comparing the expression in equation
\eqref{eq:sigma_dem} to the governing equations for a micropolar medium
\eqref{eq:goveqmicro}, we see that the conservation of torque is equivalent to the
conservation of angular momentum. Indeed, we get from \eqref{eq:Fnm} and
\eqref{eq:Fns} that
\begingroup
\def\summ{{\tiny\sum_{m=1}^{N_c}}}
\begin{equation*}
  \summ F^m \times X^m = \summ d_m^2 k_s((r - \theta)\times I^m)\times I^m = \summ d_m^2 k_s((r-\theta) - (r-\theta)\cdot I^m I^m)
\end{equation*}
so that for a square lattice ($I^1 = -I^2 = (1,0,0)$, $I^3 = -I^4 = (0,1,0)$ and $I^5
= -I^6 = (0,0,1)$), we get
\begin{equation*}
  \summ F^m \times X^m = 4d^2k_s (r - \theta).
\end{equation*}
\endgroup
and the requirement that the torque is zero yields $r - \theta=0$, which
corresponds to the conservation of angular momentum equation
\eqref{eq:angmomconsmicro}. This also highlights the need for introducing rotational
degrees of freedom for the DEM method if shear forces are used. If not, one gets the
non physical effects that rigid rotations introduce forces.  Notice that the method
which here is referred as DEM is a specific version of a lattice model where the
edges of a simplex grid are used to calculate force and the normal force is
independent of the rotation of the particles. The last statement could be understood
as neglecting rolling resistance.

We also notice that the introduction of angles has been made in finite element
literature for membrane problems \cite{Allm1984,Cook1986,Hughes1989}. In this
context, it is called the "drilling degree of freedom", see \cite{Felipp2003} for
review. The motivation has been to remove the singularity of the stiffness matrix and
the angular degree of freedom as a stiffening effect on the structure. In fact in
\cite{Hughes1995} the value for the free parameter associated with the non symmetric
part in the variation principle is recommended to be $\mu$, the shear modulo. The
degrees of freedom are completely the same as in DEM.

\subsection{Modified discrete element method MDEM}

The motivation for the introduction of the MDEM method is twofold. First, in DEM, the
relation between macro parameters and the parameters is not simple. Secondly, given a
configuration of particles, it is not possible to reproduce all the parameters
associated with isotropic materials as discussed in the paragraph above. The same
type of restrictions also holds for hexahedral and square grids, see
\cite{Suiker20011563,Pavlov20066194}. In \cite{Mott2013}, thermodynamical
considerations are used to show that, for isotropic materials, the value of the
Poisson's ratio should satisfy $\nu>1/5$. In this perspective, the restriction
$\nu<1/4$ established above for DEM appears very restrictive. The ability to vary the
mechanical properties even for this configuration introduces non central forces
between the particles, in this context called shear forces. As discussed above this
can only be done if extra local rotation variables are introduced. This has two
disadvantages, first it is more complicated, and secondly the final system is
equivalent to a micopolar medium and not a purely elastic medium. Restricting oneself
to central forces may therefore be in some cases preferable but one should remain
aware that such assumption comes with very strong restriction on the material
parameters. In \cite{Hehl2002}, the authors consider the Cauchy relations which are
known to be necessary for an elastic material where only central forces are present,
each atom or molecule is a center of symmetry and the interaction forces are well
approximated by an harmonic potential. They show that, for an isotropic material, the
Cauchy relations imply that $\nu=1/3$. This very strong restriction makes it
difficult to consider models only based on central forces.

\begin{figure}
   \includegraphics[width = \textwidth]{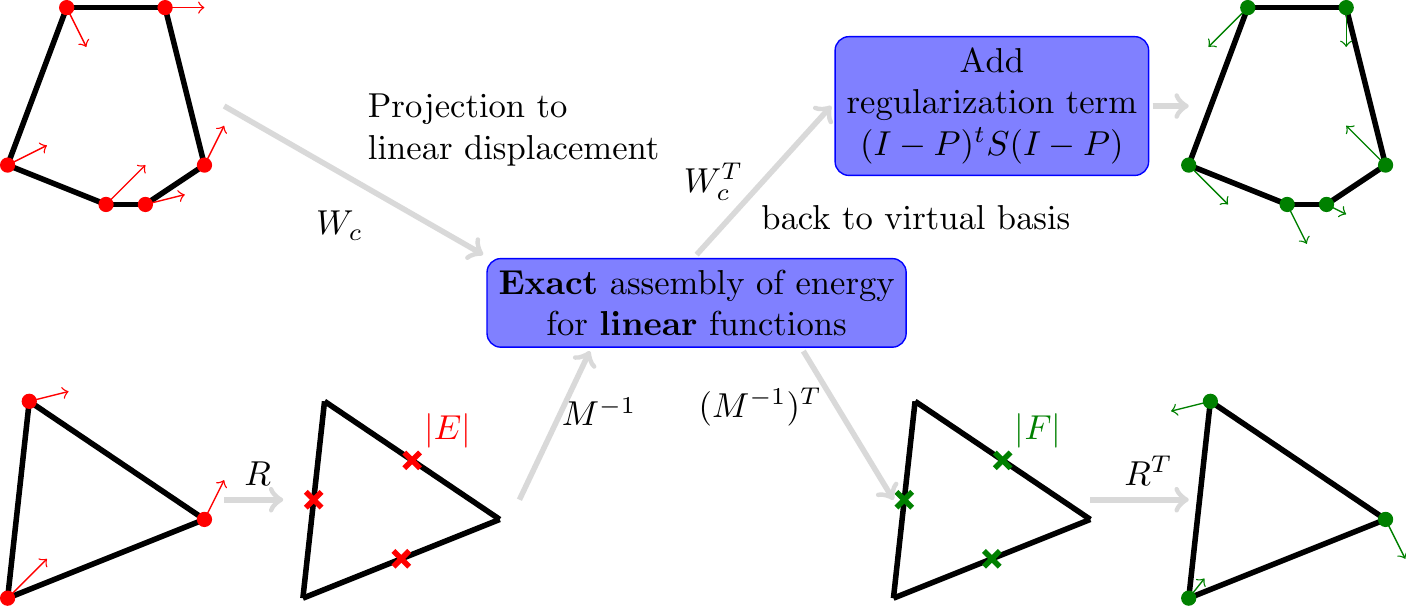}
   \caption{This show the difference and similarities of the MDEM and the VEM framework.} \label{fig:vem_dem}
\end{figure}

The basic idea of MDEM is to use an interaction region, instead of looking at the
forces on each particle as a result of interaction with neighboring particles like in
DEM, see Figure \ref{fig:vem_dem}. Then, the force at a particle is given by the sum
of the forces computed at the particle for each interactive region the particle
belongs to. In the finite element setting, the interaction region corresponds to an
element and a particle to a node. The calculation of the forces is equivalent to the
case of linear finite element. The original derivation \cite{Alassi2008,Alassi2010},
was based on explicit representation of the geometry and calculation of forces. Here
we will base our derivation on the variational form of linear elasticity. To simplify
the derivation we will use the fact that for simplex grids there exists a one to one
mapping between non rigid-body linear deformations and the length of the edges. By
non rigid-body linear deformations, we mean the quotient space of the space of linear
deformations with the space of rigid-body deformations (translation and
rotations). Such space is in bijection with symmetric matrices, the strain
tensors. Using the notation of \cite{Alassi2010} but with all tensors represented in
the Kelvin notation where the tensor inner product reduces to normal
inner-product. For simplices one can relate the non zero strain states to edge length
$U$ by
\begin{equation}
 U = M \epsi
\end{equation}
Note that, to simplify the expressions, we use different notations in the previous
section where $U$ was denoted by $\Delta U$. We write the energy of the element as
\begin{equation}
  E_{mdem} = \frac12 U^{\trans} K U = \frac12 \epsi M^{\trans} K M \epsi,
\end{equation}
where $K$ is a symmetric definite positive matrix to be determined. The tensor $K$,
which we will call in the paper the \textit{MDEM stiffness tensor}, depends on the
material parameter. This fulfills the requirement of linear elasticity that rigid
motion does not contribute to the energy. The normal forces can be calculate as the
generalized forces associated with the variable U, that is
\begin{equation}
  \label{eq:expFMDEM}
  F = \frac{\partial}{\partial U} E_{mdem} = K U.
\end{equation}
From \eqref{eq:expFMDEM}, we can see that assuming that only central forces are
present and the shear forces are negligible is equivalent to the requirement that $K$
is diagonal. Using the analogue definition of stress where we exploit the kelvin
notation
\begin{equation}
  \sigma = \frac{\partial}{\partial \epsi} \left(\frac{1}{V} E_{mdem}\right) =  \frac{1}{V} M^{\trans} K M \epsi = \frac{1}{V} M^{\trans} U.
\end{equation}
If we consider the energy of the same system for a linear elastic media assuming
constant stress and strain, which is the case for linear elements on simplex grids,
the result is
\begin{equation}
  \label{eq:Efem}
  E_{fem} = V \sigma:\epsi = V \epsi^{\trans} D \epsi
\end{equation}
where $D$ is the representation of the forth order stiffness tensor $C$ in Kelvin
notation. Note also that $epsi$ in \eqref{eq:Efem} is meant either as a tensor (in
the first equality) or as a vector written in Kelvin notations (in the second). For
the sake of simplicity, we will continue to do the same abuse of notations in the
following. We see that one reproduces the energy of a linear elastic media if
\begin{equation}
  \label{eq:relCK}
  D = M^{\trans} K M,
\end{equation}
which gives that
\begin{equation}
 K = (M^{-1})^{\trans} D M^{-1}.
\end{equation}
The difference between the matrix $K$ used in DEM and the matrix needed to reproduce
linear elasticity used in MDEM is that the latter case normally is a full
matrix. Since DEM methods solve Newtons's equation with a dissipation term it will
minimize this energy. The same is the case of standard Galerkin discretization of
linear FEM on simplices, which by construction have the same energy functional as
MDEM. Consequently, the only difference, if no fracture mechanism is present will be
the method for computing the solution to the whole system of equations. The DEM
methods rewrite the equations in the form of Newton laws with an artificial damping
term and let time evolve to converge to the solution, see \cite{Cundall1979}. For
FEM, the linear equations are usually solved directly.

The advantage of using the MDEM formulation compared to FEM is that it offers the
flexibility to choose independently on each element if a force should be computed
using linear elasticity or if a more traditional DEM calculation should be used.

The ability to associate the edge lengths to the non rigid body motions is only
possible for simplex grids. An other important aspect to this derivation is that the
degrees of freedom uniquely define all linear motions and no others. The importance
of the last part will be more evident after comparison with the VEM method.

\subsection{The Virtual Element Method}

In contrast to FEM, the Virtual element method seeks to provide consistency up to the
right polynomial order of the equation in the physical space. This is done by
approximating the bilinear form only using the definition of the degrees of freedom,
as described below.  The FEM framework on the other hand defines the assembly on
reference elements, using a set of specific basis functions. This however has
disadvantages for general grids where the mappings may be ill defined or
complicated. VEM avoids this problem by only working in physical space using virtual
elements and not computing the Galerkin approximation of the bilinear form exactly.
This comes with a freedom in the definition of the method and a cost in accuracy
measured in term of the energy norm.

As the classical finite element method, the VE method starts from the linear
elasticity equations written in the weak form of Equation
\ref{eq:lin_elast_cont},
\begin{equation}
  \label{eq:weakform}
  \int_{\Omega}\epsi( \vec{v}):C \epsi (\vec{u})  \,dx = \int_{\Omega} \vec{v}\cdot\vec{f} \,dx  \quad \text{for all} \quad \vec{v}.
\end{equation}
We have also introduced the symmetric gradient $\epsi$ given by
\begin{equation*}
  \epsi(u) = (\grad + \grad^{T} )\vec{u},
\end{equation*}
for any displacement $\vec{u}$. The fundamental idea in the VE method is to compute
on each element an approximation $a_K^h$ of the bilinear form
\begin{equation}
  \label{eq:defak}
  a_K(\vec{u}, \vec{v}) =  \int_{K} \epsi(\vec{u}) : C \epsi(\vec{v})\,dx,
\end{equation}
that, in addition of being symmetric, positive definite and coercive with respect to
the non rigid-body motions, it is also exact for linear functions. The correspondence
between MDEM and VEM we study here holds only for a first-order VEM method. When
higher order methods are used, the exactness must hold for polynomials of a given
degree where the degree determines the order of the method. These methods were first
introduced as mimetic finite element methods but later developed further under the
name of virtual element methods (see \cite{da2014mimetic} for discussions). The
degrees of freedom are chosen as in the standard finite element methods to ensure the
continuity at the boundaries and an element-wise assembly of the bilinear forms
$a_K^h$. We have followed the implementation described in \cite{gain2014}. In a
first-order VE method, the projection operator $\Pcal$ into the space of linear
displacement with respect to the energy norm has to be computed locally for each
cell. The VE approach ensures that the projection operator can be computed exactly
for each basis element. The projection operator is defined with respect to the metric
induced by the bilinear form $a_K$. The projection is self-adjoint so that we have
the following Pythagoras identity,
\begin{equation}
  \label{eq:pythagoras}
  a_K(\vec{u}, \vec{v}) = a_K(\Pcal\vec{u}, \Pcal\vec{v}) + a_K((\id - \Pcal)\vec{u}, (\id - \Pcal)\vec{v})
\end{equation}
for all displacement field $\vec{u}$ and $\vec{v}$ (in order to keep this
introduction simple, we do not state the requirements on regularity which is needed
for the displacement fields). In \cite{gain2014}, an explicit expression for $\Pcal$
is given so that we do not even have to compute the projection. Indeed, we have
$\Pcal = \Pcal_{R} + \Pcal_{C}$ where $\Pcal_{R}$ is the projection on the space $R$
of translations and pure rotations and $\Pcal_{C}$ the projection on the space $C$ of
linear strain displacement. The spaces $R$ and $C$ are defined as
\begin{align*}
  R &= \left\{ \vec{a} + B(\vec{x} - \bar{\vec{x}})\ |\ \vec{a}\in\Real^3, B\in\Real^{3\times 3},\ B^{\trans} = -B \right\},\\
  C &= \left\{ B(\vec{x} - \bar{\vec{x}})\ |\ B\in\Real^{3\times 3},\ B^{\trans} = B \right\}.
\end{align*}
Then, the discrete bilinear form $a_K^h$ is defined as
\begin{equation}
  \label{eq:disenervem}
  a_K^h(\vec{u}, \vec{v}) = a_K(\Pcal\vec{u}, \Pcal\vec{v}) + s_K((\id - \Pcal)\vec{u},(\id - \Pcal)\vec{v})
\end{equation}
where $s_K$ is a symmetric positive matrix which is chosen such that $a_K^h$ remains
coercive. Note the similarities between \eqref{eq:disenervem} and
\eqref{eq:pythagoras}. Since $\Pcal_{R}$ and $\Pcal_{C}$ are orthogonal and
$\Pcal_{R}$ maps into the null space of $a_K$ (rotations do not produce any change in
the energy), we have that the first term on the right-hand side of
\eqref{eq:pythagoras} and \eqref{eq:disenervem} can be simplified to
\begin{equation*}
  a_K(\Pcal\vec{u}, \Pcal\vec{v}) = a_K(\Pcal_C\vec{u}, \Pcal_C\vec{v}).
\end{equation*}
The expression \eqref{eq:disenervem} immediately guarantees the consistency of
the method, as we get from \eqref{eq:disenervem} that, for linear displacements,
the discrete energy coincides with the exact energy. Since the projection
operator can be computed exactly for all elements in the basis - and in
particular for the \textit{virtual} basis elements for which we do not have
explicit expressions - the local matrix can be written only in terms of the
degrees of freedom of the method. In our case the degrees of freedom of the
method are the value of displacement at the node. Let us denote
$\vec{\varphi}_i$ a basis for these degrees of freedom. The matrix
$(A_{K})_{i,j}=a_K^h(\vec{\varphi}_i, \vec{\varphi}_j)$ is given by
\begin{equation}
  \label{eq:assembvem}
  A_K =  |K| \ \ W_C^{\trans} D W_C + (\id-\Pcal)^{\trans}  S_K (\id-\Pcal).
\end{equation}
In \eqref{eq:assembvem}, $W_C$ is the projection operator from the values of node
displacements to the space of constant shear strain and $S_K$, which corresponds to a
discretization of $S_K$ in \eqref{eq:disenervem}, is a symmetric positive matrix
which guarantees the positivity of $A_{K}$. There is a large amount of freedom in the
choice of $S_K$ but it has to scale correctly.  We choose the same $S_K$ as in
\cite{gain2014}. The matrix $D$ in \eqref{eq:assembvem} corresponds to the tensor $C$
rewritten in Kelvin notations so that, in three dimensions, we have
\begin{equation*}
  D_{ij} = \epsi_i : C \epsi_j,\quad \text{ for }i,j = 1, \ldots, 6.
\end{equation*}
Finally, the matrices $A_K$ are used to assemble the global matrix $A$ corresponding
to $a^h$. In this paper, we use the implementation available as open source through
the Matlab Reservoir Simulation Toolbox (MRST) \cite{MRST:2016}. The approach of
splitting the calculation of the energy in terms of a consistent part block on one
side and a higher order block one the other side was also used in the \textit{free
  formulation} of finite elements \cite{Bergan1988}. In this case the motivation was
to find an alternative element formulation,

\subsection{Correspondences between VEM and MDEM}
\label{subsec:connec}
For simplex grids the regularization term $S_K$ in the expression for the local
stiffness matrix in Equation \ref{eq:assembvem} is zero because in this case the
projection operator is equal to the identity.  If we introduce the operator from the
degrees of freedom for the element to the edge expansions $R$, we can compare the two
expression for the local energy,
\begin{subequations}
  \label{eq:ecomp}
  \begin{equation}
    \label{eq:ecomp1}
    E_{vem} = u^{\trans} |K| W_C^{\trans} D W_C u
  \end{equation}
  and
  \begin{equation}
    \label{eq:ecomp2}
    E_{mdem} = u^{\trans} R^{\trans}  K R u =  u^{\trans} |K|  R^{\trans} (M^{-1})^{\trans} D M^{-1} R^{\trans} u.
  \end{equation}
\end{subequations}
One easily identifies the operators $W_C$ and $M^{-1} R^{\trans}$ as the projection
operator$\Pcal_c$ to the non-rigid body motions represented in Kelvin notation type
of symmetric strain. The degrees of freedom span exactly the space of linear
displacement and do not excite any higher order modes with nonzero energy. An
illustration of the different concepts is given in Figure \ref{fig:vem_dem}. We point
that both DEM and VEM calculate the basic stiffness matrix in real space, contrary to
most FEM methods which do this on the reference element.  When dealing with simplex
grids the advantage of using the DEM method within an explicit solving strategy
(often called distinct element method) is that the calculation of the edge length
extensions $U$ can be calculated for each edge, and only the matrices $M^{-1}$ and
the Cauchy stiffness tensor $C$ are needed locally. These matrices only operate on
the small space of non rigid motion with dimension $((d(d+1))/2)$ while the operator
$W_c$ works on the all the deformations which have dimension $((d+1)d)$. The edge
length can thus be seen as an efficient compact representation of the non rigid
motions, which holds only on simplices. As we have seen, both MDEM and VEM can be
derived from the calculation of the energy in each element. For the linear elastic
part it is not necessary to introduce extra angular degrees of freedom. However, this
may be needed for certain DEM methods. In this case, we refer to the use of drilling
elements in combination with the use of the free formulation of FEM
\cite{Felipp2003,Bergan1990}, which, as discussed earlier, shares some fundamental
ideas with VEM, such as \textit{energy orthogonality} (which corresponds
\eqref{eq:pythagoras}), and \textit{rc-modes exactness} (whic) and the freedom in choosing
the stabilization term.

\subsection{Fluid mechanics coupling}

We introduce a coupling with a fluid flow through the Biot's equations
\cite{biot1941general}. The Biot's equations are given by
\begin{subequations}
    \label{eq:poroelast}
    \begin{align}
      \label{eq:poroelast1}
      \dive \sigma + \grad p&= \vec{f},\\
      \label{eq:poroelast2}
      \fracpar{}{t}(S_c p + \alpha\dive \vec{u}) + \dive(\frac{K}{\mu_v}\grad p)  &= 0,
    \end{align}
\end{subequations}
where $S_c p + \alpha\dive \vec{u}$ denotes the fluid content. The fluid content
depends on the storativity $S_c$, the fluid pressure $p$ and on the rock volume
change given by $\dive\vec{u}$ which is weighted by the Biot-Willis constant
$\alpha$. In \eqref{eq:poroelast}, $K$ denotes the permeability and $\mu_v$ the fluid
viscosity. For flow, and in particular if multiphase behaviors are considered, the
most successful methods have been based on finite volume methods.  The basic time
discretization using the two point flux method or multi point flux methods \cite{
  aavatsmark2002introduction} can be written as
\begin{equation}
  S_c\frac{p^{n+1} - p^{n}}{\Delta t}  - \ddiv{f} \left[\frac{K}{\mu_v} \dgrad{p}\left[p^{n+1}\right]\right] = Q.
\end{equation}
Here $\dgrad{p}$ is a discrete gradient operator from cell pressures to face,
$\ddiv{f}$ is the corresponding discrete divergence acting on face fluxes. The source
term $Q$ represents the injection of fluids, see \cite{krogstad2015mrst} for more
details on those discrete operators.

Given an implicit time discretization the coupling term in the Biot case requires a
discrete divergence operator $\ddiv{d}$ for displacement field. Note that this
discrete operator can be implemented exactly for first order VEM, see
\cite{Raynaud:ecmor-xv} for more details. The semi-discrete equations are
\begin{equation}
  \begin{array}{c c c c c}
    A_{s} \quad \vec{u}^{n+1} &-& \alpha\ddiv{d}^{T} \left[p^{n+1}\right] &=&-F\\
    \alpha \ddiv{d} \left[\vec{u}^{n+1}\right] &+& S_c p^{n+1}-\Delta t \ddiv{f} \left[\frac{K}{\mu_v} \dgrad{p} \left[p^{n+1}\right]\right] &=& \alpha \ddiv{d} \left[\vec{u}^{n}\right] + S_c p^{n}+ Q.
  \end{array}
\end{equation}
Here $A_s$ system matrix of the mechanical system, $\ddiv{d}$ is the divergence
operator acting on the nodal displacement and gives a volume expansion of a cell and
$\alpha$ is the Biot parameter depending on the ratio between the rock and fluid
compressibility. In the context of MDEM when the simulation of fracturing is the main
purpose, we normally approximate only the volume expansion term in the transport
equation for the fractured cells, where the expansion is also the largest. Except for
this term an explicit update of pressure is used. This approximation also avoids
problems due to small permeabilities which can cause numerical locking and artificial
oscillations in the fluid pressure, see \cite{haga2012causes}.

\subsection{Fracturing criteria}

In the MDEM method, before an element is fractured, it behaves as in FEM and the MDEM
stiffness tensor $K$ is obtained from the Cauchy stiffness tensor through the
relation established in \eqref{eq:relCK}. Depending on the physical situation a
fracturing criteria based on stress is used, for example Mohr-Coulomb. In the
examples in this paper we will use the simple tensile failure criteria, namely
\begin{equation}
 \max(\sigma) > \sigma_{tens} \label{eq:frac_criteria},
\end{equation}
where $\sigma_{tens}$ is the tensile stress. After failure, we use a central force
model, where the forces are calculated individually for each edge as
\begin{equation}
 F = K_d(\Delta U) \Delta U
 \label{eq:frac_forces}
\end{equation}
where $K_d(\Delta U)$ denotes the diagonal matrix such that
\begin{equation}
   \label{eq:defKd}
  K_{d}(\Delta U)_{ii} =
  \begin{cases}
    K_{ii} &\text{if } \Delta U_i<0\\
    0&\text{otherwise}
 \end{cases}
\end{equation}
If a fracture is closing, then the effective force will in this case be as for DEM
using only central forces. As for all methods trying to simulate fracturing, the
critical point is how to avoid grid dependent fracturing, due to the singularity of
the stress field near the fracture front. In this work however, the main aim is to
see how the far-field solution can be simulated using general grid, independently of
the fracturing modeling.

\subsection{Solution method}
The solution method in MDEM is chosen to be similar to the one used in DEM.  That is
explicit time integration of Newton's laws. To get fast convergence to the physical
stationary state, the local artificial damping term that can be found
\cite{Cundall1979} is often preferable. This is not a physical damping mechanism, but
it avoids large differences in local time steps restrictions. The advantage of this
approach is that it is less sensitive to global changes than the traditional FEM
approach which solves directly the stationary state by solving the linearized
equations. This is particularly important when discontinuous changes of the forces
due to changes in the medium is present. For MDEM this is the case for the situation
of initial fracturing, equation \eqref{eq:frac_criteria} or in contact properties for
fracture cells, equation \eqref{eq:frac_forces}. The result in all cases is that the
forces are discontinuous with respect to the degrees of freedom. Explicit methods
have been shown to have advantages for such problems even if the main dynamics is
globally elliptic, because the non-linearities in the problem impose stronger
time-step size requirement that those needed for the explicit integration of the
elliptic part. As the damping criteria depends on the concept of total nodal forces,
it can also be used on the nodes connected with VEM type of force calculations. No
other modification apart from the force calculations are needed.
\begin{figure}
  \def\insertimage{
    \begin{minipage}[t]{0.3\textwidth}
      \begin{center}
        \Cval\\
        \includegraphics[width = \textwidth]{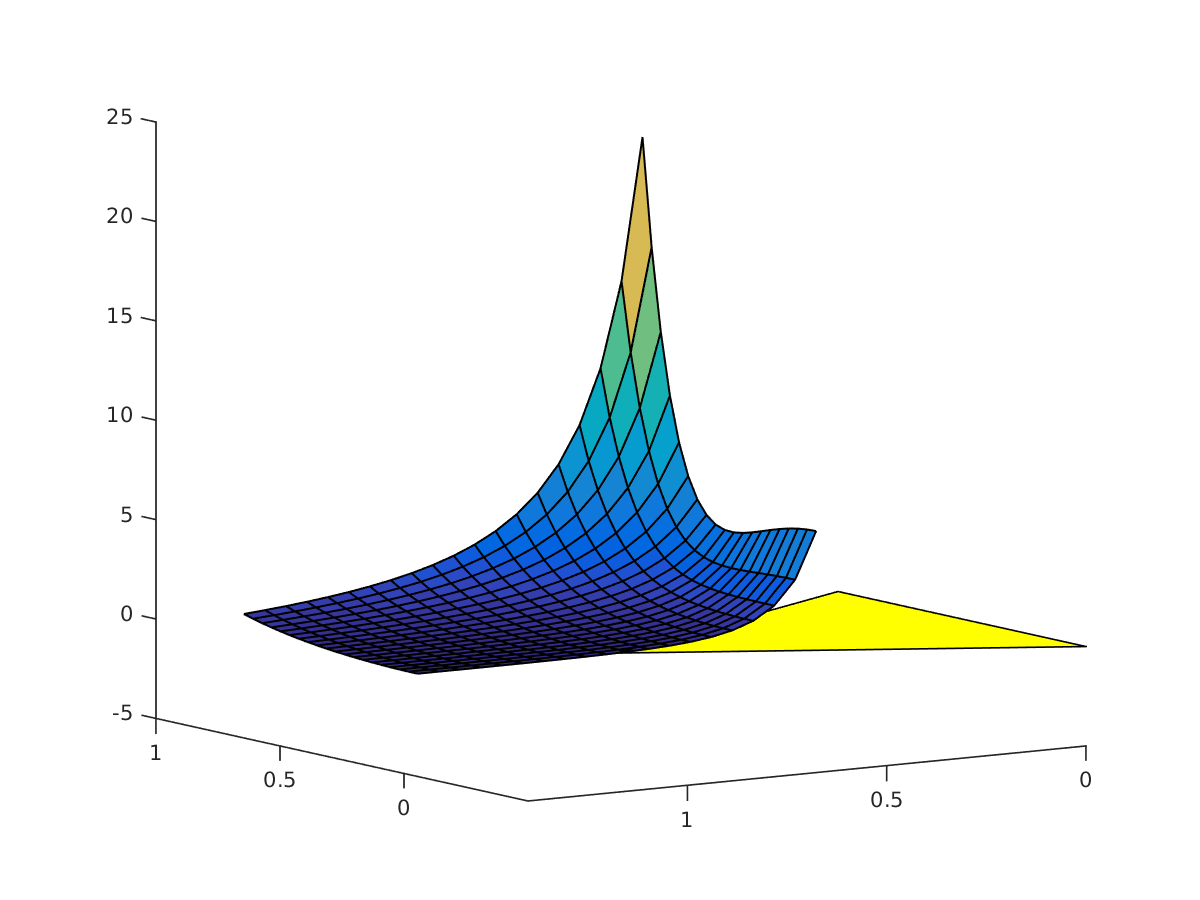}
      \end{center}
    \end{minipage}
  }
  \begin{tabular}{c c c}
    \def\imagefile{CDep_11}\def\Cval{$T_{11}$}\insertimage &
    \def\imagefile{CDep_22}\def\Cval{$T_{22}$}\insertimage &
    \def\imagefile{CDep_33}\def\Cval{$T_{33}$}\insertimage \\
    \def\imagefile{CDep_23}\def\Cval{$T_{23}$}\insertimage &
    \def\imagefile{CDep_13}\def\Cval{$T_{13}$}\insertimage &
    \def\imagefile{CDep_12}\def\Cval{$T_{12}$}\insertimage
  \end{tabular}
  \caption{The figure shows the difference of the mechanical parameters between the
    DEM fracture model \eqref{eq:frac_forces} and linear elasticity, in the case of
    compression in all edges. Here, $T = C-\Kcal^{-1}(\Kcal(C)_d)$ and the six
    components of $T$ are plotted.}
  \label{fig:C_DEM_effective_grid}
\end{figure}

\begin{figure}
\begin{tabular}{c c}
 \includegraphics[width = 0.45\textwidth]{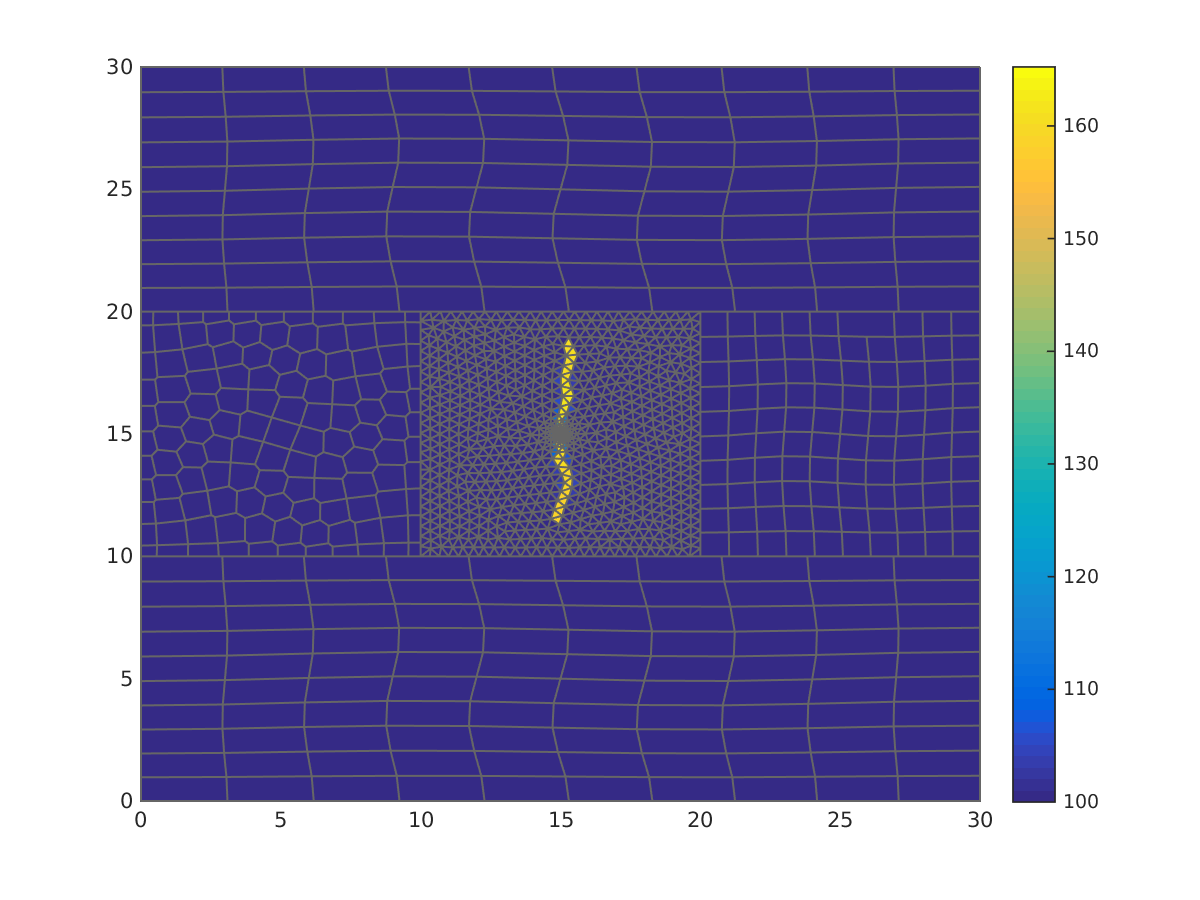} &
 \includegraphics[width = 0.45\textwidth]{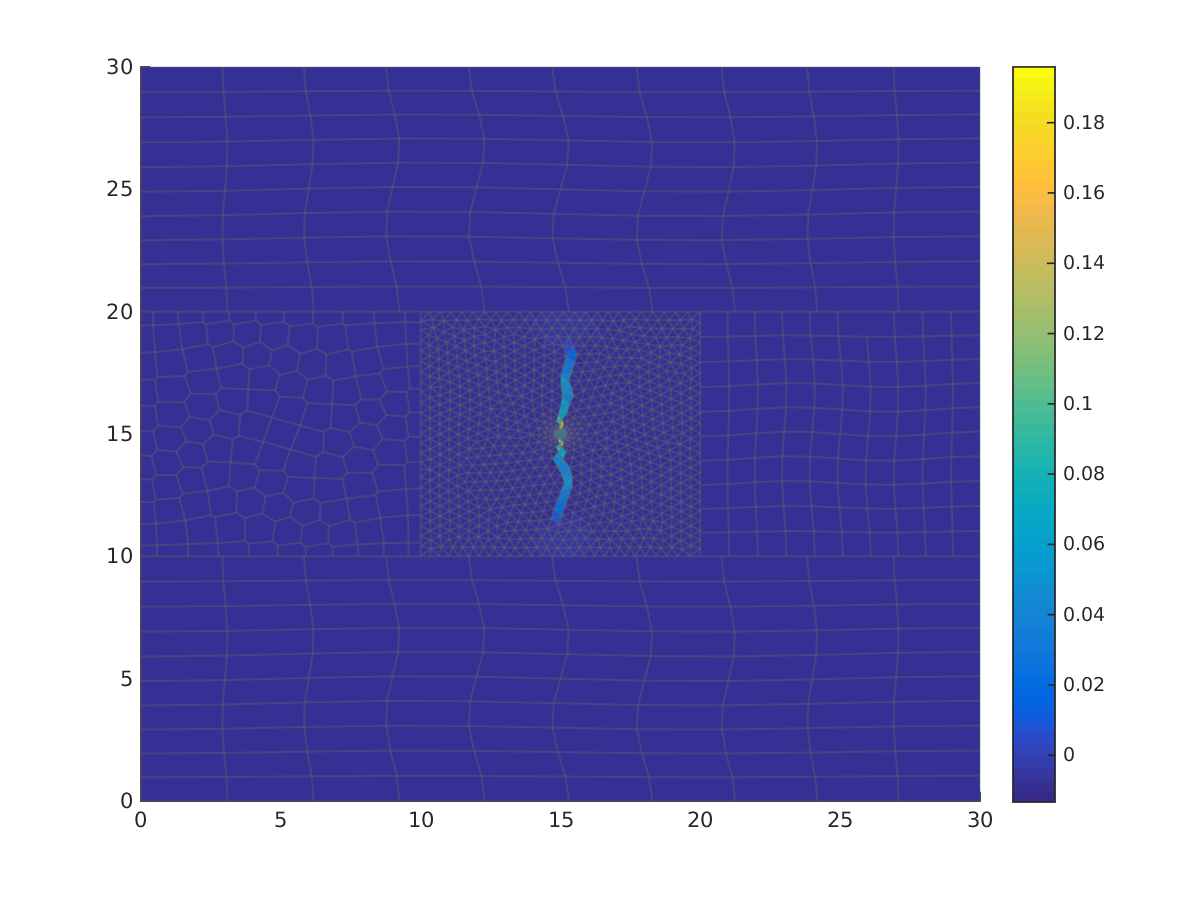} \\
 \includegraphics[width = 0.45\textwidth]{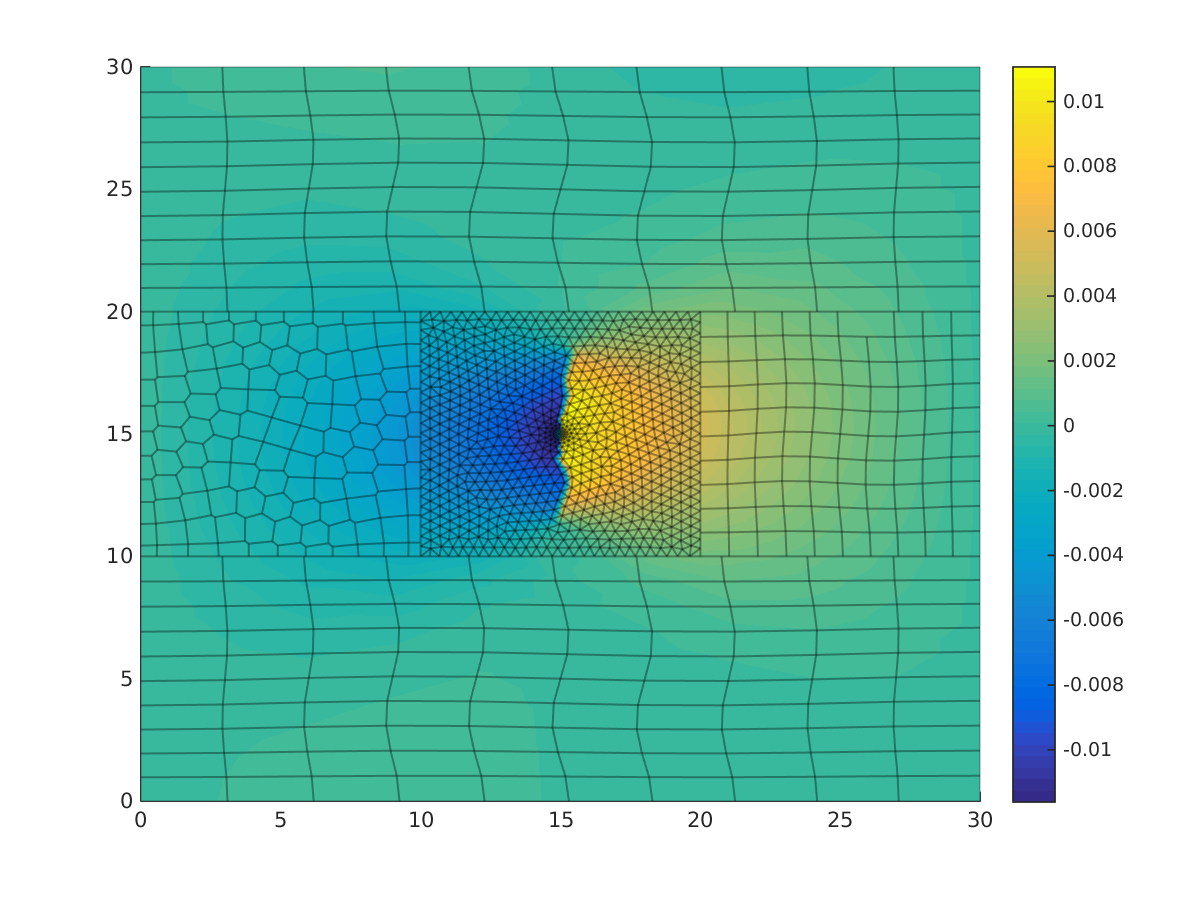} &
 \includegraphics[width = 0.45\textwidth]{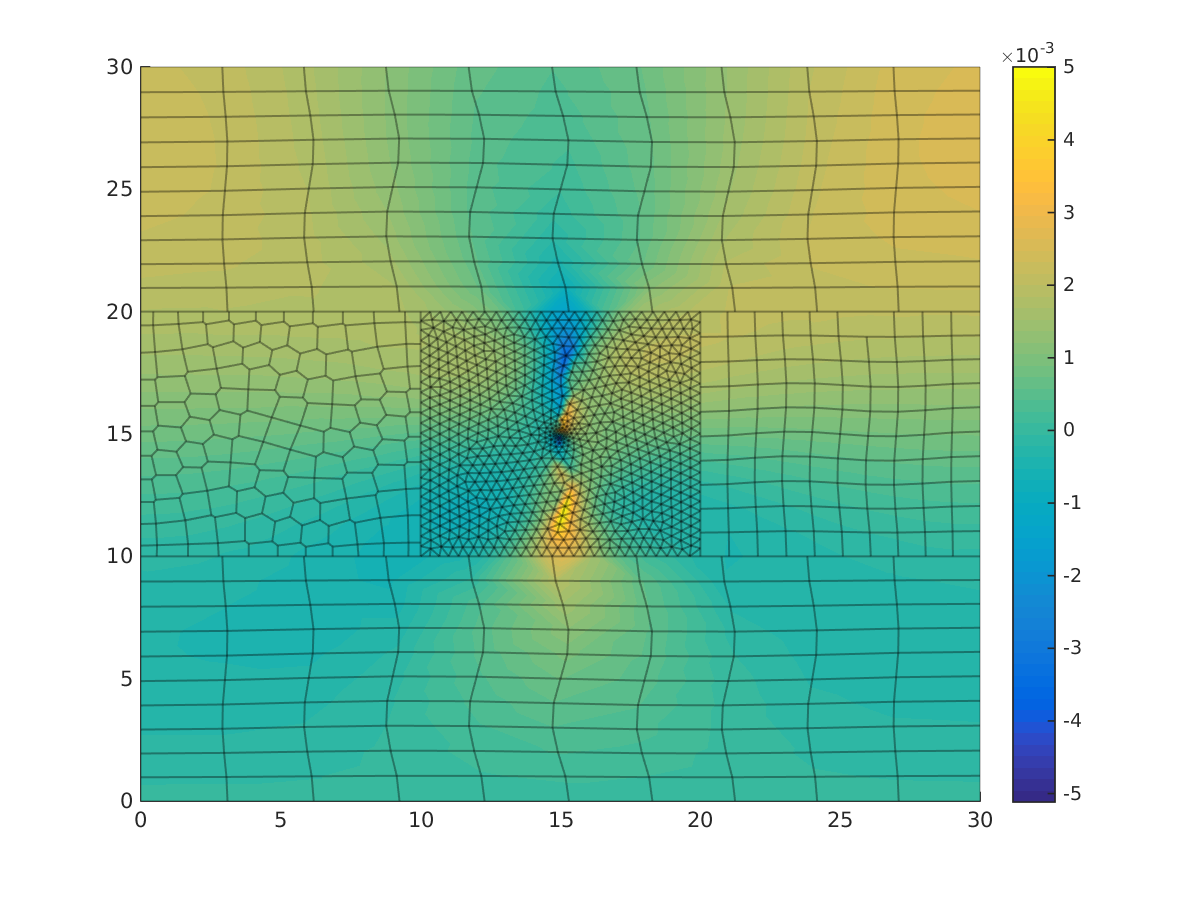} \\
 \includegraphics[width = 0.45\textwidth]{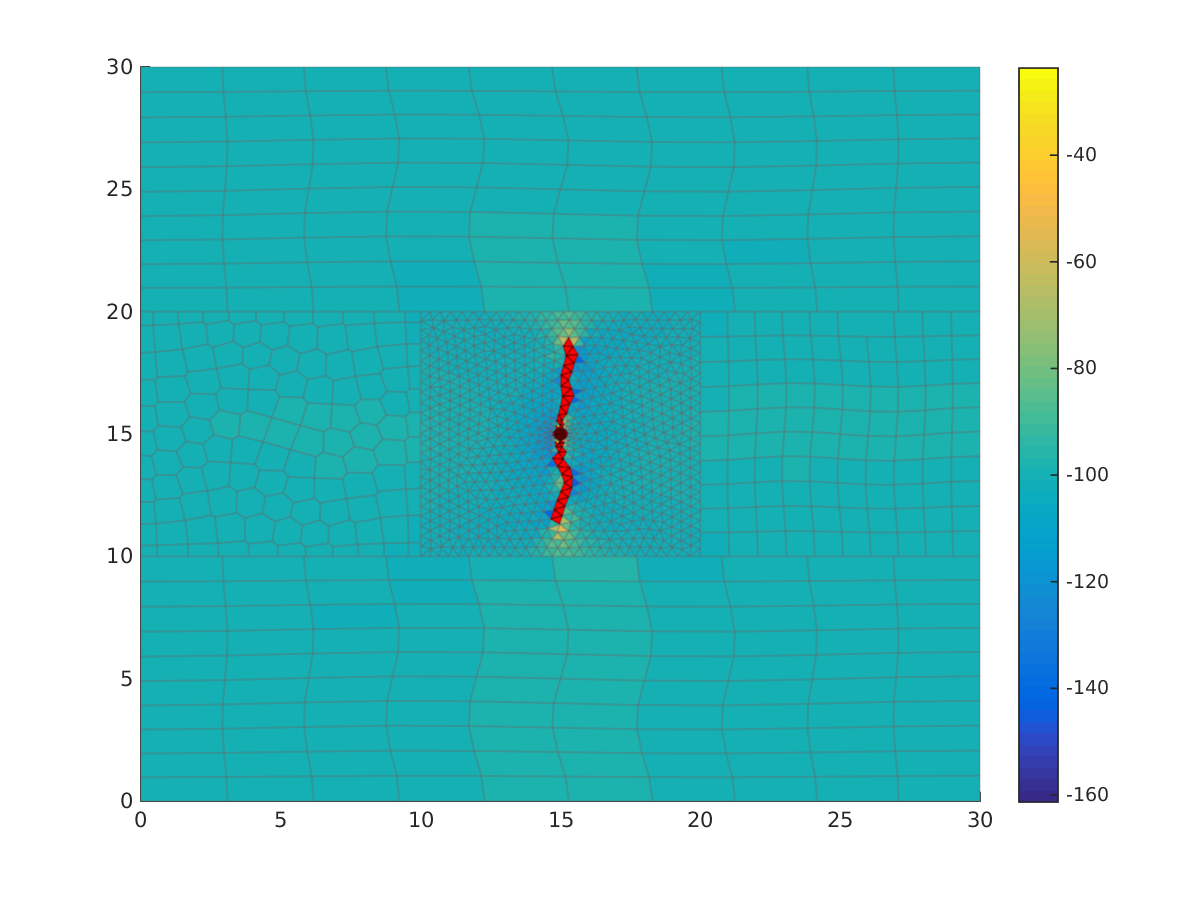} &
 \includegraphics[width = 0.45\textwidth]{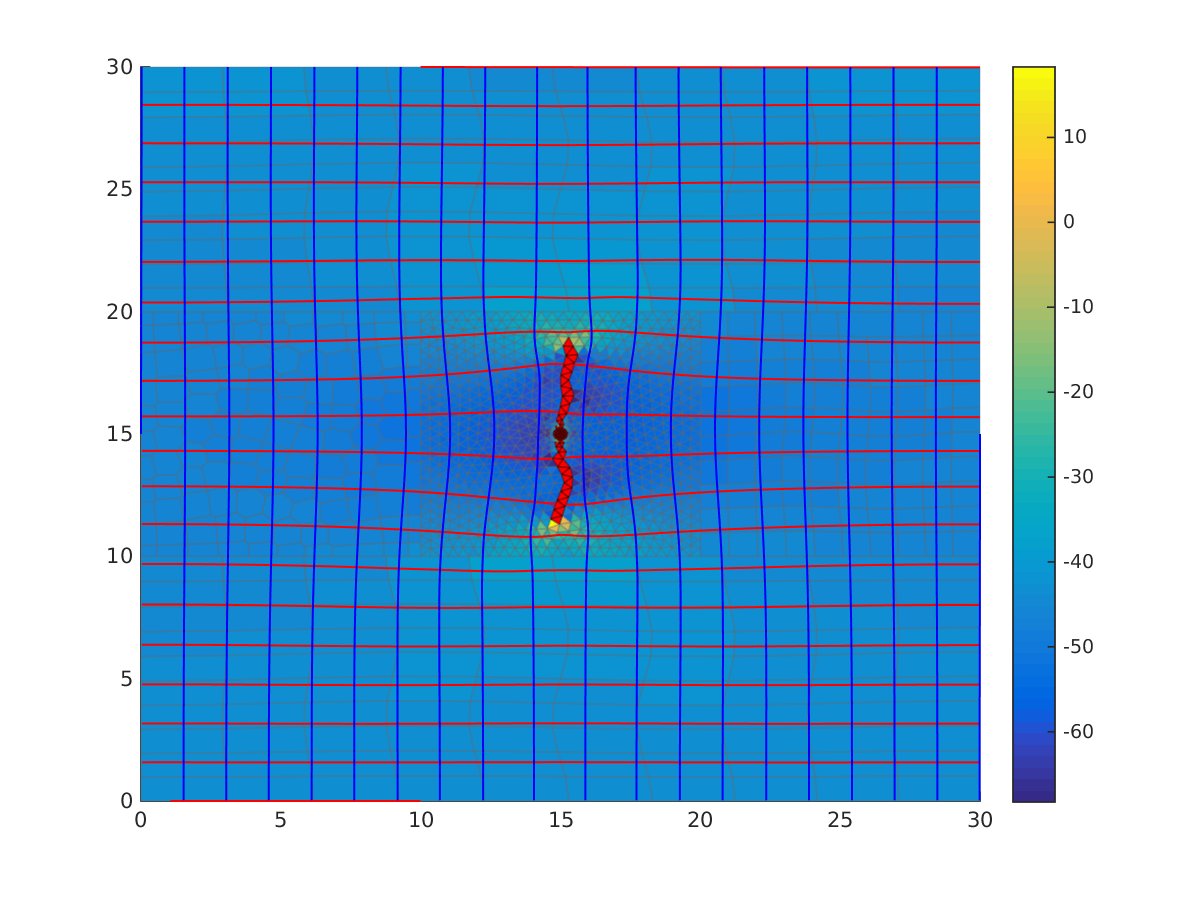}
\end{tabular}
 \caption{The figure shows a fracture growing in the direction perpendicular to
the maximum stress. Pressure in and divergence of the solution is given in upper
left and right respectively. In the middle the displacement in the x direction,
left, and y  is plotted. At the bottom the figure show the minimal
stress left and the maximum stress right. The cells in red correspond do
fractured cells. The blue and the red and blue lines is show the direction of
the principle axis for maximum and minimum stress respectively. Both pressure and stress is given in $10^5 Pa$.} \label{fig:vem_dem_fracture}
\end{figure}

\begin{figure}
\begin{tabular}{c c}
 \includegraphics[width = 0.45\textwidth]{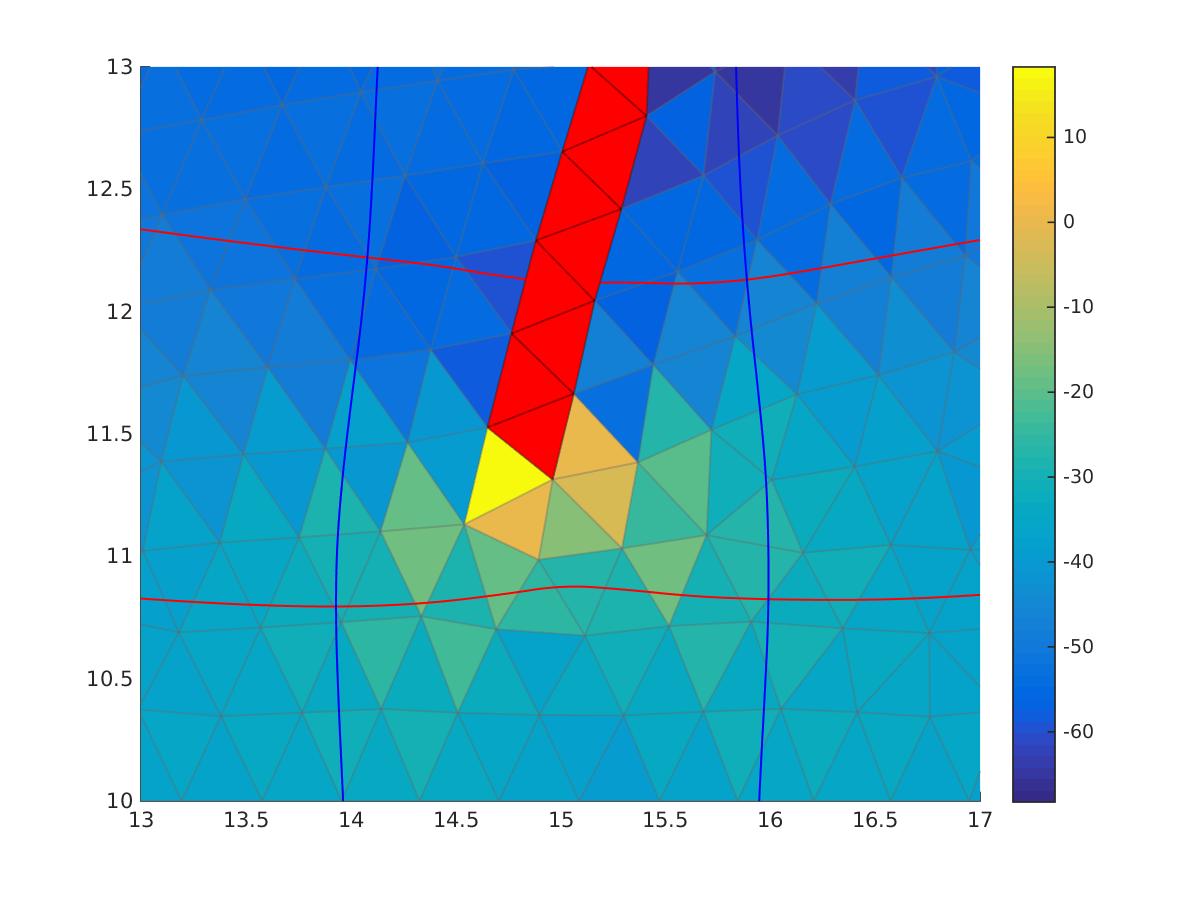} &
 \includegraphics[width = 0.45\textwidth]{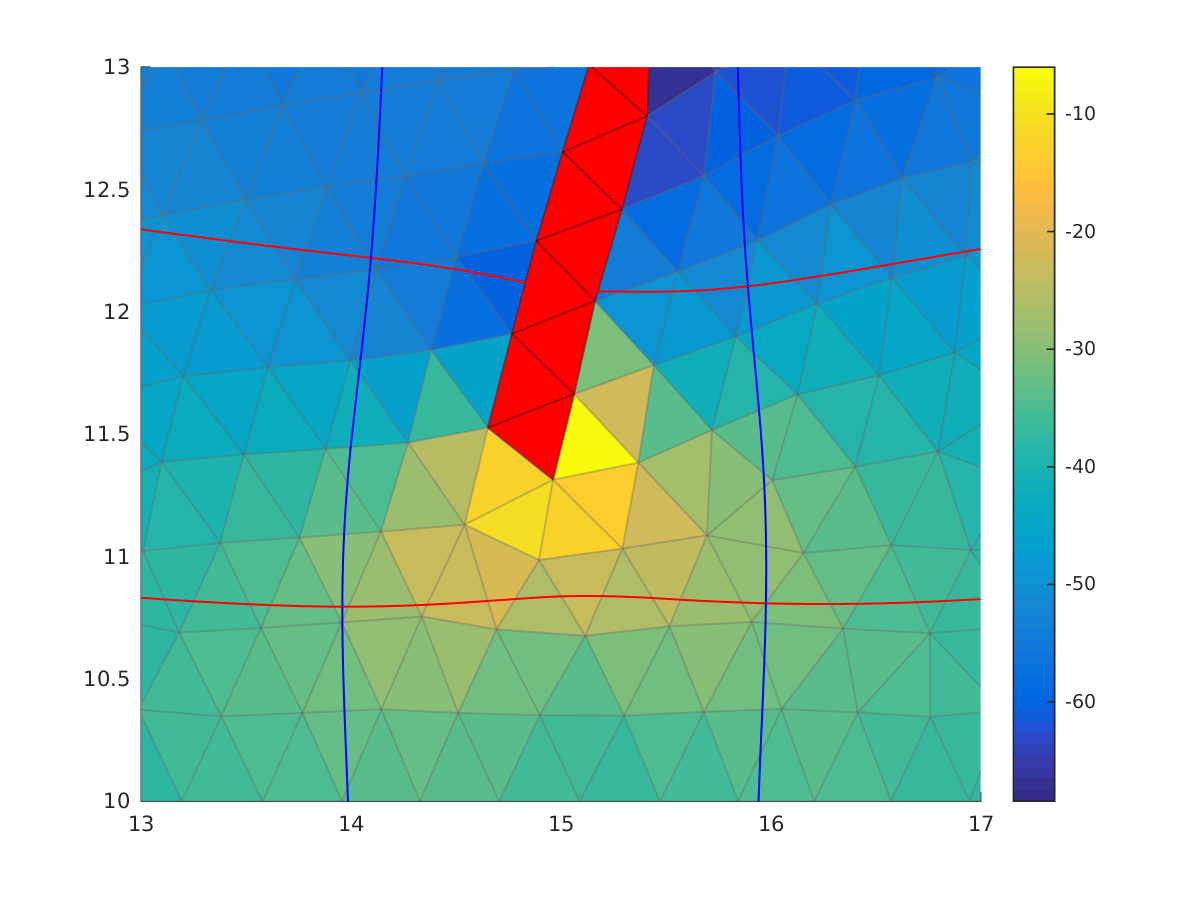}
\end{tabular}
 \caption{The figure show an large picture of the maximum stress shown in lower
right corner of figure \ref{fig:vem_dem_fracture}. The left figure are using the
stress form the linear elements while the right are from the stress using patch
recovery.}\label{fig:vem_dem_fracture_small}
\end{figure}

\section{Examples}

We demonstrate the features of the presented framework with two examples. First we
show how the effective parameters of linear elasticity in simple DEM with only
normal forces depend on the particular choice of the grid cells. Second, we use
VEM, MDEM and DEM on a general a polyhedral grid to demonstrate how this can be
combined within a uniform framework.

When a fracture has occurred in a cell, but the whole system evolves in such a way
that the fracture closes again, then we should have forces normal to the fracture
faces and, depending on the fracture model, forces along the fracture. Here, we
choose to model this by an effective stiffness tensor. Indeed, we keep using the DEM
model (and solver) after the fracture closes, meaning that the materiel parameters
for the cell are given by the diagonal MDEM stiffness tensor $K_d$ as defined in
\eqref{eq:defKd}. From Section \ref{subsec:connec}, we know that it also corresponds
to a unique Cauchy stiffness tensor. Let us study the effect of such choice and
measure the difference between the original and this \textit{post-fracturing}
stiffness tensor. If we denote by $\Kcal$, the one-to-one transformation from the
MDEM stiffness tensor $K$ to the Cauchy stiffness tensor $C$, we compute, for a given
$C$, the difference between $C$ and $\Kcal^{-1}(\Kcal(C)_d)$. We consider a
equilateral triangle and an isotropic material with Young's modulus $E = 1$ and
Poison ratio $\nu = 1/4$. For this value of $\nu$ and this shape, the matrix $K$ is
diagonal, so that $C=\Kcal^{-1}(\Kcal(C)_d)$. This reference triangle is plotted in
yellow in Figure \ref{fig:C_DEM_effective_grid}. We keep the same Cauchy stiffness
tensor but modify the shape of the triangle by translating one of the corners. For
each configuration that is obtained, we get a different post-fracturing MDEM
stiffness tensor given by $\Kcal(C)_d$ and we plot the six component of the tensor
$T$ defined as the difference $T = C-\Kcal^{-1}(\Kcal(C)_d)$.  We notice that the
changes in $C_{2,3}$ and $C_{1,3}$ is zero when on the line $x = 0$ in the
figure. This show that in this case as expected the effective model has biaxial
symmetry.  We also notice that there is quite strong changes in the effective
parameters even for relatively small changes in the triangles. It should also be
noted that a break of the edges along the x-axis, which in the MDEM fracture model
result in putting one of the corresponding diagonal element to zero only changes the
value of $C_{1,1}$. This is because this only acts in the x direction.

We use MRST \cite{MRST:2016,Lie12:comg} to generate the unstructured grid presented
in Figure \ref{fig:vem_dem_fracture} and set up an example which combines the use of
VEM for the general cell shapes and the use of MDEM for the triangular cells that can
easily be switched to a DEM model when a fracture is created. The total grid size is
$\SI{30}{\meter} \times \SI{30}{\meter}$. In the middle within a diameter
$\SI{0.5}{\meter}$ we have placed cells associated with a well. The permeability used
was $\SI{10}{\nano\darcy}$, porosity of $0.3$, the compressibility of the fluid is
similar to water, $\SI{1e-10}{\per\pascal}$, and it is injected fluid corresponding
to the pore volume of all well cells in an hour. The solution is shown after $40$
minutes. The initial condition was given by the mechanical solution with a force of
$\SI{1e7}{\pascal}$ at the top and rolling all other places. The initial condition
for pressure is constant pressure equal to $\SI{1e7}{\pascal}$. The mechanical
parameters are given by $E=\SI{1e9}{\pascal}$ and $\nu=0.3$. The well cells are set
to have Young's modulo $E=\SI{1e4}{\pascal}$ and finally the tensile strength is
$\SI{2e5}{\pascal}$. We observe that the fracture propagates in the direction so that
the fracture plane (or line in 2D) is aligned with the maximum stress plane (or line
in 2D). We get slight grid orientation effect since there is no way a planar fault in
the y direction can be obtain using the given triangular grid. The interface between
the grids has large steps in cell sizes and include hanging nodes, but no effects due
to these features are observed as long as the fracture does not reach the
interface. Near the tip of the fracture we observe oscillation of the stress on
cells, which is a well known problem for first order triangular elements. However the
values associate with the nodes is better approximated and patch recovery techniques
\cite{zienkiewicz92SPR} can be used to get better stress fields as seen in Figure
\ref{fig:vem_dem_fracture_small}. A note is that the dynamics of DEM or MDEM, is
associated with the sum of all forces from all elements around a node, not individual
stresses for cells.

\section{Conclusions}

In this paper we have shown how MDEM and VEM for linear elasticity share the same
basic idea of projection to the states of linear non-rigid motions, although with
different representations, length extension for MDEM and polynomial basis for
VEM. Both are equivalent to linear FEM on simplices, but the viewpoint presented here
gives a more direct way on how they relate. Since both share the same degrees of
freedom, except possibly the angular degree of freedom of MDEM/DEM, we combine these
methods with minimal implementation issues. This is used to simulate fracture growth,
where the near field regions is described by a simplex grid which is suited for DEM
and MDEM, while the general polyhedral grids is used in the far-field region. The
coupling between the grids which can contain hanging nodes and significant changes in
cell shapes and sizes, can be done without introducing large errors. We see this
method as a valuable contribution to flexible coupling of MDEM/DEM methods with
traditional reservoir modeling grids.

\section{Acknowledgments}

This publication has been produced with support from the KPN project
\textit{Controlled Fracturing for Increased Recovery}. The authors acknowledge the
following partners for their contributions: Lundin and the Research Council of Norway
(244506/E30).


\end{document}
